\documentclass[11pt]{amsart}

\usepackage{amsmath,amssymb,amsthm}
\usepackage{verbatim}

\usepackage{xurl}

\usepackage{xcolor}
\usepackage[psdextra]{hyperref}
\usepackage[capitalize]{cleveref}
\IfFormatAtLeastTF{2026-06-01}{%
  \newcommand*{\preamblesharedtheorem}[2]{%
    \newtheorem{#1}[theorem]{#2}%
  }%
}{%
  \usepackage{aliascnt}
  \newcommand*{\preamblesharedtheorem}[2]{%
    \newaliascnt{#1}{theorem}%
    \newtheorem{#1}[#1]{#2}%
    \aliascntresetthe{#1}%
  }%
}

\usepackage{enumerate}
\usepackage{enumitem}
\usepackage{physics}
\usepackage{todonotes}
\definecolor{crimson}{HTML}{C8102E}
\definecolor{stopgreen}{HTML}{00A14B}
\presetkeys{todonotes}{inline}{}
\presetkeys{todonotes}{color=yellow}{}

\usepackage[normalem]{ulem}

\usepackage[OT2,T1]{fontenc}
\usepackage{microtype}
\DeclareSymbolFont{cyrletters}{OT2}{wncyr}{m}{n}
\DeclareMathSymbol{\sha}{\mathalpha}{cyrletters}{"58}

\usepackage{tikz}
\usetikzlibrary{cd,decorations.markings}

\usepackage{graphicx}

\usepackage{outlines}
\usepackage{stmaryrd}
\usepackage{mathrsfs}
\usepackage{expl3}

\AddToHook{begindocument/before}{%
  \IfPackageLoadedTF{biblatex}{%
    \DeclareFieldFormat[online]{title}{\mkbibquote{#1\isdot}}%
  }{}%
}

\makeatletter
\DeclareFontFamily{OMX}{MnSymbolE}{}
\DeclareSymbolFont{MnLargeSymbols}{OMX}{MnSymbolE}{m}{n}
\SetSymbolFont{MnLargeSymbols}{bold}{OMX}{MnSymbolE}{b}{n}
\DeclareFontShape{OMX}{MnSymbolE}{m}{n}{
    <-6>  MnSymbolE5
   <6-7>  MnSymbolE6
   <7-8>  MnSymbolE7
   <8-9>  MnSymbolE8
   <9-10> MnSymbolE9
  <10-12> MnSymbolE10
  <12->   MnSymbolE12
}{}
\DeclareFontShape{OMX}{MnSymbolE}{b}{n}{
    <-6>  MnSymbolE-Bold5
   <6-7>  MnSymbolE-Bold6
   <7-8>  MnSymbolE-Bold7
   <8-9>  MnSymbolE-Bold8
   <9-10> MnSymbolE-Bold9
  <10-12> MnSymbolE-Bold10
  <12->   MnSymbolE-Bold12
}{}

\let\llangle\@undefined
\let\rrangle\@undefined
\DeclareMathDelimiter{\llangle}{\mathopen}{MnLargeSymbols}{'164}{MnLargeSymbols}{'164}
\DeclareMathDelimiter{\rrangle}{\mathclose}{MnLargeSymbols}{'171}{MnLargeSymbols}{'171}
\makeatother

\renewcommand{\phi}{\varphi}

\newcommand{\ovl}{\overline}

\renewcommand{\bar}{\overline}

\ExplSyntaxOn
\cs_new_protected:Npn \preamble_new_operator:nn #1#2
  { \cs_new:cpn {#1} { \operatorname{#2} } }
\cs_new_protected:Npn \preamble_set_operator:nn #1#2
  { \cs_set:cpn {#1} { \operatorname{#2} } }
\cs_new_protected:Npn \preamble_new_font_alias:nnn #1#2#3
  { \cs_new:cpn {#1#2} { #3 {#2} } }

\clist_map_inline:nn
  {
    ord,Gal,Spec,Hom,SL,SO,GL,Isom,Ind,Cl,Pic,End,Sym,Sel,Tor,codim,Br,Th,Lind,Id,id,
    PSL,Aut,PGL,pr,Rad,Proj,depth,Ext,gr,nil,red,rad,Vect,MT,Hdg,Sch,Frob,Ver,Fr,fr,
    Vol,Div,Supp,cl,Sw,Coim,Gr,Ker,Coker,coker,Cov,Mod,Lie,Stab,acl,dcl,Mult,eq,ACF,
    Kr,Ded,Av,PD,st,Ch,Var,Swan,val,rk,Der,SU,SH,Ann,Cone,sel,tp,stp,supp,tor,gp,
    disc,Cn,Jac,Vis,WC,Comp,Frac,Func,Perv,Sp,diag,Nrd,an,Hilb,Cliff,Ad,Diff,sing,
    Out,Bl,rec,Inn,Loc,Coh,Aff,colim,coim,Ob,Sub,Def,alg,nr,un,ad,Fl,ab,fl,Ric,loc,
    sm,Mor,aut,Hol,An,Arg,BGL,BU,BO,full,om,Sing,St,Lk,ch,td,per,ind
  }
  { \preamble_new_operator:nn {#1} {#1} }

\preamble_new_operator:nn {Char} {char}
\preamble_new_operator:nn {Maxspec} {MaxSpec}
\preamble_new_operator:nn {ext} {Ext}
\preamble_new_operator:nn {height} {ht}
\preamble_new_operator:nn {vol} {Vol}
\preamble_new_operator:nn {rM} {r-M}
\preamble_new_operator:nn {im} {im}
\preamble_new_operator:nn {Sh} {SH}
\preamble_new_operator:nn {cone} {Cone}
\preamble_new_operator:nn {defect} {def}
\preamble_new_operator:nn {resid} {res}
\preamble_new_operator:nn {Graff} {GrAff}

\preamble_set_operator:nn {Form} {Form}
\preamble_set_operator:nn {div} {div}
\preamble_set_operator:nn {ker} {Ker}
\preamble_set_operator:nn {Im} {Im}
\preamble_set_operator:nn {df}{def}
\ExplSyntaxOff

\newcommand{\SHom}{\mathcal{H}\kern -.5pt\mathit{om}}
\newcommand{\shom}{\mathcal{H}\kern -.5pt\mathit{om}}
\newcommand{\sext}{\mathcal{E}\kern -.5pt\mathit{xt}}

\newcommand{\xr}{\xrightarrow}

\renewcommand{\subset}{\subseteq}

\ExplSyntaxOn
\clist_map_inline:nn {A,B,C,D,E,F,G,H,I,J,K,L,M,N,O,P,Q,R,S,T,U,V,W,X,Y,Z}
  {
    \preamble_new_font_alias:nnn {b}  {#1} {\mathbf}
    \preamble_new_font_alias:nnn {bb} {#1} {\mathbb}
    \preamble_new_font_alias:nnn {s}  {#1} {\mathscr}
    \preamble_new_font_alias:nnn {mc} {#1} {\mathcal}
  }
\clist_map_inline:nn {a,b,c,d,e,f,g,h,i,j,k,l,m,n,o,p,q,r,s,t,u,v,w,x,y,z}
  { \preamble_new_font_alias:nnn {mf} {#1} {\mathfrak} }
\ExplSyntaxOff

\catcode`,\active

\catcode`\,12

\newtheorem{theorem}{Theorem}[section]

\crefname{theorem}{Theorem}{Theorems}
\Crefname{theorem}{Theorem}{Theorems}
\crefname{thm}{Theorem}{Theorems}
\Crefname{thm}{Theorem}{Theorems}

\preamblesharedtheorem{proposition}{Proposition}
\crefname{proposition}{Proposition}{Propositions}
\Crefname{proposition}{Proposition}{Propositions}

\preamblesharedtheorem{prop}{Proposition}
\crefname{prop}{Proposition}{Propositions}
\Crefname{prop}{Proposition}{Propositions}

\preamblesharedtheorem{lemma}{Lemma}
\crefname{lemma}{Lemma}{Lemmas}
\Crefname{lemma}{Lemma}{Lemmas}

\preamblesharedtheorem{lem}{Lemma}
\crefname{lem}{Lemma}{Lemmas}
\Crefname{lem}{Lemma}{Lemmas}

\newtheorem*{lemma*}{Lemma}

\preamblesharedtheorem{claim}{Claim}
\crefname{claim}{Claim}{Claims}
\Crefname{claim}{Claim}{Claims}

\preamblesharedtheorem{corollary}{Corollary}
\crefname{corollary}{Corollary}{Corollaries}
\Crefname{corollary}{Corollary}{Corollaries}

\preamblesharedtheorem{conjecture}{Conjecture}
\crefname{conjecture}{Conjecture}{Conjectures}
\Crefname{conjecture}{Conjecture}{Conjectures}

\newtheorem{innerimportantthm}{Theorem}

\crefname{innerimportantthm}{Theorem}{Theorems}
\Crefname{innerimportantthm}{Theorem}{Theorems}

\newtheorem{innerimportantcor}{Corollary}

\crefname{innerimportantcor}{Corollary}{Corollaries}
\Crefname{innerimportantcor}{Corollary}{Corollaries}

\theoremstyle{definition}
\preamblesharedtheorem{definition}{Definition}
\crefname{definition}{Definition}{Definitions}
\Crefname{definition}{Definition}{Definitions}

\preamblesharedtheorem{nondefinition}{Non-Definition}
\crefname{nondefinition}{Non-Definition}{Non-Definitions}
\Crefname{nondefinition}{Non-Definition}{Non-Definitions}

\preamblesharedtheorem{exercise}{Exercise}
\crefname{exercise}{Exercise}{Exercises}
\Crefname{exercise}{Exercise}{Exercises}

\preamblesharedtheorem{example}{Example}
\crefname{example}{Example}{Examples}
\Crefname{example}{Example}{Examples}

\preamblesharedtheorem{remark}{Remark}
\crefname{remark}{Remark}{Remarks}
\Crefname{remark}{Remark}{Remarks}

\preamblesharedtheorem{observation}{Observation}
\crefname{observation}{Observation}{Observations}
\Crefname{observation}{Observation}{Observations}

\preamblesharedtheorem{question}{Question}
\crefname{question}{Question}{Questions}
\Crefname{question}{Question}{Questions}

\preamblesharedtheorem{problem}{Problem}
\crefname{problem}{Problem}{Problems}
\Crefname{problem}{Problem}{Problems}

\numberwithin{equation}{section}

\usepackage[style=alphabetic,maxalphanames=4]{biblatex}
\begin{document}

\pagestyle{plain}

\title{$\bbC(t)$ has the strict order property}
\author{Spencer Dembner and Cole Franklin}
\date{November 9, 2026}

\begin{abstract}
Let $(E, P)$, $(E', P')$ be elliptic curves over $\bbC$ with marked $p$-torsion points; a \emph{marked isogeny} is an isogeny sending $P$ to $P'$. We show the marked isogeny relation is interpretable in $\bbC(t)$ and has strict infinite chains. It follows that $\bbC(t)$ has the strict order property, which confirms a special case of the stable fields conjecture.
\end{abstract}

\maketitle

\section{The main theorem}

One of the central questions in model theory is how its tameness conditions interact with the theory of fields. Macintyre~\cite{Macintyrecategoricalfields} showed that an infinite field $(K, +, \cdot)$ in the language of rings is \emph{$\omega$-stable} if, and only if, $K$ is algebraically closed. Cherlin and Shelah~\cite{CherlinShelahSuperstableFields} later showed the same is true if $K$ is merely \emph{superstable}. However, we currently have no characterization of \emph{stable} fields. If $K$ is separably closed then $(K, +, \cdot)$ is stable; the \emph{stable fields conjecture} predicts that the converse also holds.

The field $\bbC(t)$, which is not separably closed, should therefore be unstable. But showing this has proven difficult, and in recent years $\bbC(t)$ has become a notable test case. Johnson, Tran, Walsberg, and Ye~\cite{JTWYEtaleOpenTopology} proved the stable fields conjecture for \emph{large} fields, a class that excludes $\bbC(t)$, and suggested $\bbC(t)$ could be a counterexample to the general conjecture. There do exist many interesting sets definable in $\bbC(t)$. For example, the constant field $\bbC \subset \bbC(t)$ is definable, because a positive genus curve $X / \bbC$ has no nonconstant $\bbC(t)$-points. More strikingly, the set $\mathbf{CM} \subset \bbC$ of $j$-invariants for elliptic curves with complex multiplication, a countable dense subset of $\bbC$, is definable. But using results related to the André-Oort conjecture, Scanlon~\cite{ScanlonDecidability} showed that the structure $(\bbC, +, \cdot, \mathbf{CM})$ is stable (and decidable).

Our main theorem resolves this uncertainty.

\begin{theorem}\label{thm}
    The theory $\Th(\bbC(t))$ is unstable. In fact, it has the strict order property.
\end{theorem}

\begin{proof}
    In this section we describe a relation $R$ that witnesses the strict order property, and show it has strict infinite chains. The proof that $R$ is interpretable requires somewhat more care, and is postponed to the next section.

    Let $p > 0$ be prime, and consider pairs $(E,P)$ where $E / \bbC$ is an elliptic curve and $P \in E[p]$ is a $p$-torsion point. Define the \emph{marked isogeny relation} $R$ in the following way: We have $R((E,P), (E', P'))$ if there  is an isogeny
    \[ f \colon E \to E', \qquad f(P) = P'. \]
    We will show in the next section that $R$ is interpretable, and clearly $R$ is reflexive and transitive. Therefore, we only need to exhibit a strict infinite chain. Let $E / \bbC$ be a non-CM elliptic curve, and let $\{P, Q \}$ be a basis for $E[p]$. Write $Q = Q_0$, and choose a series of points $Q_i$ with $p Q_{i+1} = Q_i$. Let $E_i = E / \langle Q_i \rangle$ and let $P_i \in E_i[p]$ represent the image of $P$. If $i < j$ then the canonical projection map $\pi \colon E_i \to E_j$ is a marked isogeny. Suppose for contradiction that there is a marked isogeny $f$ in the other direction. Then, since $E_i$ is not CM, the composite map
    \[ E_i \xr{\pi} E_j \xr{f} E_i  \]
    is the multiplication map $[n]$ for some $n$. Since $\deg(\pi)$ is a multiple of $p$, $n$ is also a multiple of $p$. It follows that $[n] P_i = 0$, which contradicts the fact that $\pi$ and $f$ are marked isogenies. \qedhere
\end{proof}

It remains to show that the marked isogeny relation is interpretable. Scanlon \cite[Section 3]{ScanlonDecidability} showed that the \emph{unmarked} isogeny relation is interpretable and in fact much more: If $E_s, E_t$ denote the elliptic curves with $j$-invariant $s$ and $t$, one can uniformly interpret the groups
\[ \Hom_{\bbC}(E_s, E_t)  = E_t(\bbC(E_s)) / E_t(\bbC). \]
If \emph{composition} of isogenies were interpretable, then one could interpret the endomorphism ring $\End(E) \cong \bbZ$ of a non-CM elliptic curve. It would follow that $\bbC(t)$ is unstable, and indeed undecidable as well. But it is unclear whether this should be possible: In this setting isogenies appear as function field-valued points, not as composable maps. Instead, we will show that composition is definable ``at finite level:'' For fixed $N$, one can uniformly interpret the maps
\[ E[N] \to E'[N] \]
induced by isogenies. This trivially suffices to interpret the marked isogeny relation.

\begin{remark}
    For the proof of \cref{thm} it would be enough to interpret the marked isogeny relation only on $2$-torsion points. It is possible to do this directly using coordinates on the Legendre family
    \[ y^2 = x (x-1)(x - \lambda), \]
    but we will not pursue this.
\end{remark}

\subsection*{Future directions}

The other longstanding problem about $\bbC(t)$ is whether its theory is undecidable, and it appears that the methods here are not enough to settle this. Scanlon \cite[Question 5.5]{ScanlonDecidability} has also asked whether the structure $(\bbC, +, \cdot, \rm{Isog})$ generated by the bare isogeny relation is stable, and this also remains unanswered.

\subsection*{Acknowledgments} The authors thank Ravi Vakil and Tom Scanlon for helpful discussions.

\subsection*{AI Declaration}

GPT-6 Astra found the first proof of \cref{thm}. The authors developed the ``level structure'' point of view in \cref{sec:Interp}, simplified other aspects of the proof, and wrote the paper.

\section{Isogenies with level structure}\label{sec:Interp}

If $N \geq 3$ is an integer, then the \emph{modular curve} $Y(N)$ is an algebraic variety over $\bbC$. It is the moduli space for pairs
 $(E, \alpha)$ where $E$ is an elliptic curve and
\[ \alpha \colon (\bbZ / N)^2 \xrightarrow{\sim} E[N] \]
trivializes its $N$-torsion. Since $\bbC \subset \bbC(t)$ is definable, the set $Y(N)(\bbC)$ is definable. There is a universal family
\[ \mcE(N) \to Y(N) \]
of elliptic curves over $Y(N)$, whose fiber over $(E, \alpha)$ is the elliptic curve $E$, and there are canonical sections $Y(N) \to \mcE(N)$ which uniformly trivialize the torsion subgroups $E[N]$.

\begin{proposition}\label{prop:torsion}
    Let $a$ and $b$ range over $Y(N)(\bbC)$. There is a uniformly interpretable family of abelian groups $H_{a,b}$, representing $\Hom_{\bbC}(E_a, E_b)$, and the maps
    \[ E_a[N] \to E_b[N] \]
    induced by elements of $H_{a,b}$ are uniformly definable.
\end{proposition}

We will prove \cref{prop:torsion} at the end of the section. We first require the following lemma, which shows coherent sheaves give definable families over a field.

\begin{lemma}\label{lem:coherent}
    If $X$ is an algebraic variety over a field $k$, and $\mcF$ is a coherent sheaf on $X$, then the fibers
    \[ \{ \mcF(x): x \in X(k) \} \]
    are a uniformly interpretable family of vector spaces. If $\phi:\mcF\to \mcF'$ is an $\mcO_X$-linear map of coherent sheaves, then
    \[\phi(x):\mcF(x)\to \mcF'(x)\]
    are a uniformly interpretable family of $k$-linear maps. If $\mcA$ is a coherent sheaf of $\mcO_X$-algebras, then $\mcA(x)$ is a uniformly interpretable family of $k$-algebras.

\end{lemma}

\begin{proof}
    First, suppose that $X \subset \bbA^n_k$ is an affine variety. Then there exists a global presentation
    \[ \mcO_X^n \xr{f} \mcO_X^m \to \mcF \to 0. \]
    The map $f$ is represented by a matrix $A$ whose entries are polynomial functions on $x$, and the fiber $\mcF(x)$ is the cokernel of $A(x)$. This shows that the fibers form an interpretable family. If $\mcF'$ is another coherent sheaf, then it also admits a presentation
    \[ \mcO_X^{n'} \to \mcO_X^{m'} \to \mcF' \to 0, \]
    and since $X$ is affine, the map $\mcF \to \mcF'$ lifts to a map $\mcO_X^m \to \mcO_X^{m'}$. This map is polynomial, so the induced map $\mcF(x) \to \mcF'(x)$ is $k$-interpretable.

    If $X$ is a general variety, then it admits a finite open cover by affine varieties $U_i$, with $k$-definable transition functions. The previous argument shows that the fibers $\mcF(x)$ are uniformly interpretable over each $U_i$, and the definable transition functions make these interpretations compatible on overlaps. As above, a map $\mcF \to \mcF'$ has polynomial representatives on each chart, and these induce compatible maps $\mcF(x) \to \mcF'(x)$ for points $x$ in the overlap of two charts.

    The statement about the algebra structure follows from the previous cases applied to the coherent sheaves $\mcA$ and $\mcA\otimes \mcA$ and the multiplication map $\mcA\otimes_{\mcO_X}\mcA \to \mcA.$ \qedhere
\end{proof}

\begin{lemma}\label{lem:families}
    Let $X, Y$ be varieties over $\bbC$, and let
    \[ f \colon Y \to X \times \bbP^1 \]
    be a finite morphism whose fibers over points of $X$ are integral curves. The function fields
    \[ \{ \bbC(Y_x): x \in X  \} \]
    are a uniformly interpretable family of field extensions over $\bbC(t)$. Every $f$-equivariant automorphism $\sigma \in Aut(Y)$ gives a uniformly interpretable family of field automorphisms over $\bbC(t)$.
\end{lemma}

\begin{proof}
    Since $\Spec(\bbC(t))$ is the generic point of $\bbP^1$, there is a natural map
    \[  i \colon X \times_{\bbC} \Spec \bbC(t) \to X \times \bbP^1. \]
    Since $f$ is finite, $i^*f_*\mcO_Y$ is a coherent sheaf on $X\times_\bbC\Spec \bbC(t).$
    By \cref{lem:coherent}, its fibers
    \[  (i^* f_* \mcO_Y)(x), x \in X(\bbC) \]
    form a uniformly interpretable family. Since $f$ is finite, we have the following pullback diagram:

\[\begin{tikzcd}
	{\Spec \bbC(Y_x)} & {Y_x} & Y \\
	{\Spec \bbC(t)} & {\{x\}\times \bbP^1} & {X\times \bbP^1}
	\arrow[from=1-1, to=1-2]
	\arrow[from=1-1, to=2-1]
	\arrow["\lrcorner"{anchor=center, pos=0.125}, draw=none, from=1-1, to=2-2]
	\arrow[from=1-2, to=1-3]
	\arrow[from=1-2, to=2-2]
	\arrow["\lrcorner"{anchor=center, pos=0.125}, draw=none, from=1-2, to=2-3]
	\arrow[from=1-3, to=2-3]
	\arrow[from=2-1, to=2-2]
	\arrow[from=2-2, to=2-3]
\end{tikzcd}\]

But $f$ is also an affine morphism, and it follows that we have
\[ i^* f_* \mcO_Y(x) = \bbC(Y_x). \]
Similarly, since an $f$-equivariant $\sigma\in \text{Aut}(Y)$ gives an automorphism $\sigma:i^*f_*\mcO_Y\to i^*f_*\mcO_Y,$ \cref{lem:coherent} gives a uniformly interpretable family of automorphisms $\sigma(x): \bbC(Y_x)\to \bbC(Y_x)$. \qedhere

\end{proof}

\begin{proof}[Proof of \cref{prop:torsion}]
    Assume that $N>2.$ Consider the modular curve $Y(N)$, and the universal elliptic curve $\mcE(N) \to Y(N)$. Let $(P, Q)$ be the distinguished basis for $E[N]$ coming from the level structure. Then there is a canonical map
    \[ \mcE(N) \to Y(N) \times \bbP^1. \]
    Indeed, quotienting by $[-1]$ gives a map to a $\bbP^1$-bundle on $Y(N)$, which is trivialized by sending $(O,P,Q)$ to $(0, 1, \infty)$. This map is proper with zero-dimensional fibers, and therefore finite. Precomposing with multiplication by $N$ gives a map
    \[ \mcE(N) \xrightarrow{[N]} \mcE(N) \to Y(N) \times \bbP^1, \]
    equivariant for all the translations by $N$-torsion points. By \cref{lem:families}, we can uniformly interpret the fields $\bbC(E)$ in this family as extensions of $\bbC(t)$, and we can uniformly interpret the field automorphisms $\sigma_P$ arising from translation by $P \in E[N]$.  \\

    Now take $a,b\in Y(N)(\bbC),$ and let $f\in E_b(\bbC(E_a))$ be thought of as a map $f:E_a\to E_b.$ If $\sigma_P$ is the automorphism arising from translation by the point $P \in E_a[N]$, and $\ovl f$ is the unique homomorphism arising as an $E_b(\bbC)$-translate of $f$, then we have

    \[\sigma_P(f)-f=\bar f(P)\in E_b(\bbC)[N]\subset E_b(\bbC(E_a)). \]
    Letting $H_{a,b} = E_b(\bbC(E_a)) / E_b(\bbC)$, we obtain the desired maps $H_{a,b}\to \Hom(E_a[N], E_b[N]).$
    \end{proof}

\begin{corollary}
For $a,b \in Y(N)(\bbC)$, the marked isogeny relation $R((E_a, P_a),(E_b, P_b))$ is interpretable.
\end{corollary}

\printbibliography

@article {Macintyrecategoricalfields,
    AUTHOR = {Macintyre, Angus},
     TITLE = {On {$\omega \sb{1}$}-categorical theories of fields},
   JOURNAL = {Fund. Math.},
  FJOURNAL = {Polska Akademia Nauk. Fundamenta Mathematicae},
    VOLUME = {71},
      YEAR = {1971},
    NUMBER = {1},
     PAGES = {1--25},
}

@article {ScanlonDecidability,
    AUTHOR = {Scanlon, Thomas},
     TITLE = {Decidability of some complicated structures definable in {$\mathbb{C}(t)$}},
   JOURNAL = {J. Symbolic Logic},
  FJOURNAL = {The Journal of Symbolic Logic},
      YEAR = {2026},
      NOTE = {Accepted manuscript, published online 30 October 2026},
       DOI = {10.1017/jsl.2026.10279},
}

@article {CherlinShelahSuperstableFields,
    AUTHOR = {Cherlin, Gregory L. and Shelah, Saharon},
     TITLE = {Superstable fields and groups},
   JOURNAL = {Ann. Math. Logic},
  FJOURNAL = {Annals of Mathematical Logic},
    VOLUME = {18},
      YEAR = {1980},
    NUMBER = {3},
     PAGES = {227--270},
       DOI = {10.1016/0003-4843(80)90006-6},
}

@article {JTWYEtaleOpenTopology,
    AUTHOR = {Johnson, Will and Tran, Chieu-Minh and Walsberg, Erik and
              Ye, Jinhe},
     TITLE = {The \'{e}tale-open topology and the stable fields conjecture},
   JOURNAL = {J. Eur. Math. Soc.},
  FJOURNAL = {Journal of the European Mathematical Society},
    VOLUME = {26},
      YEAR = {2024},
    NUMBER = {10},
     PAGES = {4033--4070},
       DOI = {10.4171/JEMS/1345},
}

\end{document}